\documentclass[final, 4p,10pt]{elsarticle}
\usepackage{amsmath}
\usepackage{array}
\usepackage{longtable}
\newcommand{\oan}{\mbox{OA}(N,k,s,t)}

\newcommand{\oal}{\mbox{OA}(\lambda s^t,k,s,t)}
 \newcommand{\oapl}{\mbox{OAP}(k,s,t, \lambda)}
 \newcommand{\gloapl}{\Pi \mbox{(OAP}(k,s,t, \lambda))}

\newcommand{\N}{{\bf N}}

\newcommand{\1}{{\bf 1}}
\newcommand{\0}{{\bf 0}}
\newcommand{\aaa}{{\bf a}}
\newcommand{\bb}{{\bf b}}

\newcommand{\uu}{{\bf u}}

\newcommand{\x}{{\bf x}}
\newcommand{\vv}{{\bf v}}
\newcommand{\y}{{\bf y}}
\newcommand{\w}{{\bf w}}

\newcommand{\z}{{\bf z}}

\newcommand{\B}{{\bf B}}
\newcommand{\ddd}{{\bf d}}
\newcommand{\D}{{\bf D}}

\newcommand{\A}{{\bf A}}

\newcommand{\rr}{\mathbb{R}}

\usepackage{amssymb}

\usepackage{simplemargins}
\setleftmargin{2in}
\setrightmargin{2in}
\settopmargin{1.9375in}
\setbottommargin{1.9375in}

\newtheorem{thm}{Theorem}
\newtheorem{lem}{Lemma}

\newtheorem{defn}{Definition}
\newdefinition{rmk}{Remark}
\newproof{pf}{Proof}
\newproof{pot}{Proof of Theorem \ref{thm2}}
\journal{JCMCC}

\begin{document}

\begin{frontmatter}



\author{A.J. Geyer}
\ead{andrew.geyer@afit.edu}

\address{Air Force Institute of Technology/ENC, 2950
Hobson Way WPAFB, OH 45433-7765.}

\author{D.A. Bulutoglu\corref{cor1}}
\ead{dursun.bulutoglu@afit.edu}
\cortext[cor1]{Corresponding author.}
\address{Air Force Institute of Technology/ENC, 2950
Hobson Way WPAFB, OH 45433-7765.}

\author{S.J. Rosenberg}
\ead{srosenbe@uwsuper.edu}

\address{Mathematics and Computer Science Department,
 University of Wisconsin Superior,  Swenson Hall 3023,
Belknap and Catlin
P.O. Box 2000
Superior, WI 54880.}

\title{The LP Relaxation Orthogonal Array Polytope and its  Permutation Symmetries}


\begin{abstract}
Symmetry plays a fundamental role in design of experiments. In particular, symmetries of factorial designs that preserve their statistical properties
are exploited to find designs with the best statistical properties. By using a result proved by Rosenberg~\cite{Rosenberg1995},
 the concept of the LP relaxation orthogonal array polytope is developed and studied. A complete characterization of the   permutation symmetry group of this polytope is made.  Also, this characterization is verified computationally for many cases. Finally, a proof is provided. 

\end{abstract}


\end{frontmatter}


\section{Introduction}
\label{}
A factorial design $\D$ with $N$ runs, $k$ factors, and
$s \geq 2$ levels is an orthogonal array of strength $t$, where $1 \leq t
\leq k$, (denoted by OA($N,k,s,t$)) if each of the $s^t$ $t$-level
combinations appear exactly $\lambda=N/s^t$ times when $\D$ is projected
onto any $t$ factors. Orthogonal arrays are known to be universally optimal for estimating certain statistical models.
However, for a given $N,k,s,t$ combination finding an OA($N,k,s,t$) or proving that it does not exist is a notoriously difficult problem.

Two orthogonal arrays are isomorphic if one can be obtained from the other by permuting factors or runs as well as
permuting levels in a subset of factors. 
Bulutoglu and Margot~\cite{Bulutoglu2008} enumerated all non-isomorphic orthogonal arrays for many $N$, $k$, $s$, $t$ combinations. 
Their enumeration was based on finding all non-isomorphic, non-negative integer solutions to a system of equations with binary coefficients, where each variable represents the number of 
times a factor level combination occurs in the sought after orthogonal array and each solution corresponds to an orthogonal array.
To find all such solutions, they used the integer programming solver with isomorphism pruning developed by Margot~\cite{Margot2007}.
 They declared two solutions to be isomorphic if one could be obtained from the other by applying a permutation  belonging to the automorphism group of the constraint matrix. In each enumeration in~\cite{Bulutoglu2008}, the Margot~\cite{Margot2007} integer programming solver used this automorphism group.
  
For a $k$ factor, $s$ level design $\D$, Bulutoglu and Margot~\cite{Bulutoglu2008}  defined the isomorphism group of $\D$, $G_{s,k}$ to be the group of all permutations of factors as well as
permutations of levels of factors in $\D$. $G_{s,k}\cong S_s\wr S_k$, where ``$\wr$" is the wreath product and $|S_s \wr S_k|=(s!)^kk!$. 
Bulutoglu and Margot~\cite{Bulutoglu2008} showed that the $G_{s,k}$  of an $\oal$ is a subgroup of the automorphism group $G$ of the constraint matrix defining orthogonal arrays with the same parameters. They also observed that $G_{s,k}=G$ for each $\oal$ they enumerated.
In this paper, $G_{s,k}=G$ unless $t=0$ or $t=k$ is shown. Furthermore, all the  permutation symmetries in the linear programming (LP) relaxation of the orthogonal array problem are given.

The  system of equations used in~\cite{Bulutoglu2008} for finding orthogonal arrays had also been used by Rosenberg~\cite{Rosenberg1995} to prove the following lemma.
\begin{lem}\label{lem:Ros}
Define a sequence $\{a_c\}$ recursively by 

\[
a_0=\lambda, \quad a_c=\lambda-\sum_{e=0}^{c-1} a_e {k-t \choose c-e}(s-1)^{c-e} \quad \text{for } c \geq 1. 
\]
Let $\z$, $\x$ and $\y \in \{1,\ldots,s\}^k$ be row vectors with $0 \leq d(\z,\x)\leq t$ where $d(\z,\x)$ is the number of non-zero entries in $\z-\x$.
Also, let $I_{\x}=\{i \in \{1,\ldots,k\}:x_i \neq z_i\}$ and $J_{\x}=\{ \y \in  \{1,\ldots,s\}^k:y_i=x_i \ \ \forall i \in I_{\x} \}$.
 Then

\begin{equation} \label{eqn:Ros}
\begin{aligned}
&N_{\x}=a_{t-d(\z,\x)} +(-1)^{t-d(\z,\x)+1}\sum_{\y\in J_{\x} \atop d(\z,\y)>t} { d(\z,\y)-d(\z,\x)-1 \choose t-d(\z,\x)} N_{\y},\\ 
&N_{\y}\geq0, \text{ for $\y$ such that $d(\z,\y)>t$},
\end{aligned}
\end{equation}
where $N_{\x}$, $N_{\y}$ are the number of times factor level combinations $\x$, $\y$  appear in a conjectured $\oal$.
\end{lem}

In section \ref{sec:pol}, we use Lemma \ref{lem:Ros} to define the LP relaxation orthogonal array polytope of an $\oal$.
We then show that an $\oal$ exists if and only if its LP relaxation orthogonal array polytope denoted by $\oapl$ contains integer vectors where each integer vector represents an $\oal$. Furthermore, we prove that our formulation defining $\oapl$ does not contain any distinct redundant inequalities which renders each distinct inequality a facet.

In section \ref{sec:sym},  the group $\gloapl$ consisting of all the  permutation symmetries of $\oapl$    is found and  completely characterized.
 This characterization is numerically verified for many cases.
Finally, this characterization and $G_{s,k}=G$ unless $t=k$ or $t=0$ are proven.

\section{The LP relaxation orthogonal array polytope} \label{sec:pol}
The following theorem follows immediately from Lemma \ref{lem:Ros} by taking $\z=\1$ and observing that $N_{\x}\geq0$.
\begin{thm}\label{thm:orthogonalarray}
Let  $\{a_c\}_{c=0}^t$ and $I_{\x}$ be as in Lemma \ref{lem:Ros}. Let $0 \leq d(\1,\x)\leq t$ 
and
\begin{equation}\label{eqn:Rosin}
\begin{aligned}
a_{t-d(\1,\x)} +(-1)^{t-d(\1,\x)+1}\sum_{\y\in J_\x \atop d(\1,\y)>t} { d(\1,\y)-d(\1,\x)-1 \choose t-d(\1,\x)} N_{\y}\geq 0,&\\
N_{\y}\geq 0,&
\end{aligned}
\end{equation}
where $J_{\x}=\{ \y \in  \{1,\ldots,s\}^k:y_i=x_i \ \ \forall i \in I_{\x} \}$. 
Then an $\oal$ exists if and only if there exist integers $N_{\y}$ satisfying the inequalities (\ref{eqn:Rosin}). 
\end{thm}
Theorem \ref{thm:orthogonalarray} converts the Bulutoglu and Margot~\cite{Bulutoglu2008}  integer linear programming (ILP) problem with equalities to an ILP problem with  
inequalities only. This is done by deleting the set of basic variables $N_{\x}$ with $d(\x,\1)\leq t$ after implementing Gaussian elimination. Deletion of these variables is only possible because the coefficients of all $N_{\x}$ are $1$ and the coefficients of $N_{\y}$ with $d(\y,\1)> t$ are all integers at the end of  Gaussian elimination.
The substance of Lemma \ref{lem:Ros} in~\cite{Rosenberg1995} is that it correctly identifies which set of basic variables is an integer combination of the remaining variables.
This enables us to delete these variables from the Bulutoglu and Margot~\cite{Bulutoglu2008} ILP. In Section \ref{sec:sym},  it is shown that 
the size of the symmetry group that can be exploited by the solver in~\cite{Margot2007} is $k!((s-1)!)^k$ if variables are deleted and $k!(s!)^k$ otherwise.
A speed comparison of these two formulations under Margot~\cite{Margot2007} solver where their corresponding symmetry groups were exploited  is made in Table 1.
For each $\oan$ enumerated, the second and third columns of Table 1  report the number of enumerated $\oan$ by exploiting the groups of size $k!(s!)^k$
and $k!((s-1)!)^k$  before and after deleting variables respectively. Likewise, fourth and fifth columns report solution times before and after variables are deleted.  
Even though the new formulation has fewer variables, computational experiments summarized in Table 1  suggest that it should not be preferred over the original formulation.  
 It is evident from Table 1 that exploiting the larger symmetry group more than overcomes the additional computational burden of having a larger number of variables. 
In fact, the computational savings appear to grow exponentially with the number of variables. On the other hand, the cases OA$(64,7,2,4)$ and OA$(24,11,2,3)$ do buck this trend.

For a general ILP with only equality constraints over non-negative integer vectors, deleting variables by using Gaussian elimination may not always be possible either because there is no set of basic variables that are integer combinations of the remaining variables plus some integer or it may be very difficult to identify such a set  of basic variables.
For example, for the ILP in~\cite{Bulutoglu2008} when $k=8$, $s=2$, and $t=3$, we estimated the proportion of such sets of basic variables to all sets of basic variables to be $.5\%$. This estimate was calculated by repeating the following procedure $1000$ times. First randomly permute the columns of the constraint matrix, augment the resulting matrix with its right hand side, then  row reduce it  to its reduced row echelon form and record if the output has only integer entries.

Let $m=\sum_{i=(t+1)}^k {k \choose i} (s-1)^i$. The following definition arises naturally from Theorem \ref{thm:orthogonalarray}.
\begin{defn}\label{def:def}
The set of $\N^{y} \in \rr^m$ with $k \geq d(\1,\y) \geq t+1$  that  satisfy the system of inequalities (\ref{eqn:Rosin}) is called the LP relaxation polytope of $\oal$ denoted by $\oapl$.
\end{defn}
Note that OAP$(k,s,t, \lambda)$ is a polytope not an unbounded polyhedron because it can be embedded inside the hypercube $[0,\lambda]^m$.
\begin{thm} \label{thm:full}
 The $\oapl$  is full dimensional for all $\lambda$.
\end{thm}
\begin{pf}
Let $m$ be as above and $\1_m$ be the $m \times 1$ vector of all $1$s. We prove this result by showing that $\lambda/s^{k-t}\1_{m}$  is an interior point of the $\oapl$. Rosenberg~\cite{Rosenberg1995} showed that the system of equations in Lemma \ref{lem:Ros} is equivalent to
\begin{equation}\label{eqn:mar}
\sum_{x \in [s]^k\atop x_i=a_i \forall i \in I}N_{\x}=\lambda \quad \text{for all $t$-subsets $I$ of $[k]$ and all $\aaa \in [s]^k$,} 
\end{equation}
where $[s]^k=\{1,\ldots,s\}^k$ and $[k]=\{1,\ldots,k\}$.
 It is clear that $\lambda/s^{k-t}\1_{s^k}$ solves the system of equations (\ref{eqn:mar}) in $({\mathbb{R}^+})^{s^k}$.
Then $N_{\x}=\lambda/s^{k-t}$ and $N_{\y}=\lambda/s^{k-t}$ for $d(\x,\1)\leq t$ and  $d(\y,\1)> t$ respectively solves  the  Lemma \ref{lem:Ros} system of constraints (\ref{eqn:Ros}). Now, clearly  $\lambda/s^{k-t}\1_{m}$ satisfies all the inequalities defining the $\oapl$ strictly. Hence it is an interior point.
\end{pf}

Next, we prove that none of the constraints in Theorem \ref{thm:orthogonalarray} is redundant unless $k=t+1$ and $s=2$.

\begin{thm} \label{conj:facet}
Each one of the distinct $s^k$ inequalities in Theorem \ref{thm:orthogonalarray} defining the $\oapl$ is a facet and no facet is repeated unless 
  $k=t+1$ and $s=2$.
\end{thm}
\begin{pf}
At least one of the inequalities in (\ref{eqn:Rosin}) is a facet since otherwise $\oapl$ would be an unbounded polyhedron. Then there exists  $\N^{y} \in \mathbb{R}^m$ satisfying all but the facet defining inequality in (\ref{eqn:Rosin}). Let $N^{f(\uu_0)}_{\uu_0}<0$, $N^{f(\uu_0)}_{\uu_0} \in \mathbb{R}$ be the left hand side in constraints (\ref{eqn:Ros}) corresponding to the facet defining inequality in (\ref{eqn:Rosin}), where
\begin{equation*}
f(\uu)= \begin{cases}
  x& \mbox{if } d(\uu,\1)\leq t,\\
 y&  \mbox{otherwise.}  
\end{cases}
\end{equation*}
Hence there exist  vectors  $\N^{y} \in \mathbb{R}^m$ and $\N^{x} \in \mathbb{R}^{s^k-m}$ that satisfy the equations in constraints (\ref{eqn:Ros})
such that  $N^y_{\y} \geq 0$ and $N^x_{\x} \geq 0$ for $\x\neq \uu_0$ and $\y\neq \uu_0$, where  $N^{f(\uu_0)}_{\uu_0} < 0$ . 
The group $G_{s,k} \cong S_s\wr S_k$ sends vectors in $\mathbb{R}^{s^k}$ that satisfy the equations in
constraints  (\ref{eqn:Ros}) to vectors that satisfy the same equations. Furthermore $G_{s,k}$ acts transitively on the variables of 
constraints  (\ref{eqn:Ros}). Hence, for each $\uu_0\in 
[s]^k$, there exists a solution with $N^{f(\uu_0)}_{\uu_0}<0$ and  $N^{f(\w)}_{\w} \geq 0$ for $\w\neq \uu_0$, $\w \in [s]^k$. Then there exists $\N^{y} \in \mathbb{R}^m$ satisfying all but one facet defining inequality in (\ref{eqn:Rosin}) whose left hand side is $N^{f(\uu_0)}_{\uu_0}<0$ in constraints (\ref{eqn:Ros}) for arbitrary $\uu_0 \in [s]^k$.
Hence, there are no distinct redundant inequalities in (\ref{eqn:Rosin})  and each distinct inequality is a facet.

 Observe that unless $k=t+1$ and $s=2$,  $$\{\y \in J_{\x_1} : d(\1,\y)>t\}\neq \{\y \in J_{\x_2} : d(\1,\y)>t\}$$ whenever $\x_1\neq 
\x_2$. Hence no facet is repeated unless 
  $k=t+1$ and $s=2$.
For the degenerate case  $k=t+1$ and $s=2$, there is only one variable $N^y_{(2,2,\ldots,2)}$ and we get $N^y_{(2,2,\ldots,2)} \geq 0$ and $-N^y_{(2,2,\ldots,2)} \geq -\lambda$ each repeated  $2^{k-1}$ times. 
\end{pf}
\begin{rmk}
While it is true that  $ S_s\wr S_k$ acts as a group of symmetries on   $\oapl$, this action is no longer as a group of linear transformations  (as is the case for the full system 
of equations (\ref{eqn:mar})), but rather as a group of affine transformations. In particular $ S_s\wr S_k$ does not permute the variables of $\oapl$, but rather, it acts 
by permuting the half-spaces defined by inequalities (\ref{eqn:Rosin}). Furthermore, this action is transitive.
\end{rmk}

We verified Theorem \ref{conj:facet} when $\lambda=1$ for each of $s=2$, $4\leq k\leq 13$, $2 \leq t\leq k-2$,
 $s=3$, $3\leq k\leq 8$, $2 \leq t\leq k-1$, and  $s=4$, $3\leq k\leq 6$, $2 \leq t\leq k-1$ cases.
Our verification was based on finding interior points  on each facet of the $\oapl$ as follows: 
Let $\B \x \leq \ddd$ be the system of inequalities in  Theorem \ref{thm:orthogonalarray} defining the $\oapl$.
Let $\B^i \x \leq \ddd^i$ be the same system after the $i$'th inequality is deleted. Also, let $(\bb_i)^T$ be the $i$'th row of $\B$  and $F_i$ be the hyperplane defined by the equality 
$(\bb_i)^T\x=d_i$.
To find interior points  on $F_i \cap \oapl$, each face of the $\oapl$, find feasible solutions to the following linear program
\begin{eqnarray*}\label{eqn:ILP}
&\min \quad \1^T\x \nonumber \\
&{\rm such \ that} \quad (\bb^i)^T\x = d_i \\
& \quad \B^i\x \leq \ddd^i-\frac{1}{1000}\1. \nonumber
\end{eqnarray*}

\section{Permutation symmetries of the LP relaxation orthogonal array polytope}\label{sec:sym}
We first define the  permutation symmetries of a polytope. 
\begin{defn}\label{defn:intsym}
Let $P$ be a full dimensional polytope in $ {\mathbb{ R}}^m$. A
permutation of coordinates of ${\mathbb {R}}^m$  that also sends $P$ onto itself is called 
a {\em permutation symmetry} of $P$. The set of all such transformations forms a group called the {\em  permutation symmetry group  ($\Pi(P)$)} of $P$.
\end{defn}

The variables $N_{\y}$ in Theorem \ref{thm:orthogonalarray} are indexed by all factor level combinations $\y \in [s]^k$ with $d(\y,\1)>t$ and $\gloapl$ permutes those variables.
The following theorem explicitly describes a nontrivial subgroup of $\gloapl$.

\begin{thm}\label{thm:Hkst}
Let 
\begin{equation*}
 H_{k,s,t} \cong
\begin{cases}
S_{s^k-1}&  \text{if $t=0$,}\\
 S_{(s-1)} \wr  S_k & \text{if $0<t<k$ and ($k>t+1$ or $s>2$)},\\
I & \text{otherwise,}
\end{cases}
\end{equation*}
  and $I$ be the identity group. Also, let $Y=\{\y \in [s]^k: d(\y,\1)>t\}$. 
Then, when $H_{k,s,t}$ is not defined to be $I$, it naturally embeds as the group of permutations that preserve $Y$,   
and $H_{k,s,t}\subseteq\gloapl$.
\end{thm}
\begin{pf}
It is easy to see that each element of $H_{k,s,t}$  maps Definition \ref{def:def}  defining constraints of the $\oapl$ to each other. Then, $\vv \in \oapl$ $ \Rightarrow h(\vv) \in \oapl$
for all $h \in H_{k,s,t}$. Hence $H_{k,s,t}$ maps the $OAP(k,s,$ $t,\lambda)$ into itself. On the other hand, for a given $\vv \in \oapl$, we have $h(h^{-1}(\vv))=\vv$ as
$h^{-1}(\vv) \in \oapl$. This implies that $h(OAP(k,s,t,$ $\lambda))=\oapl$ for each $h$. Hence  $H_{k,s,t} \subseteq \gloapl$. 
\end{pf}
\begin{rmk}\label{rem:remark} For $k=t+1$ and $s=2$ there is only one variable $N_{\y}$ with $d(\y,\1)>t$, hence  $\gloapl=I$.
For the case $t=0$, there is one constraint on the $s^k-1$ variables. The  coefficients of this constraint are all $-1$s, hence 
$H_{k,s,0}\cong S_{s^k-1}$. If $t=k$ there are no variables in inequalities (\ref{eqn:Rosin}).
\end{rmk}






Next, we develop tools for calculating $\gloapl$.  
As noted in~\cite{Margot2010}, the set of  all permutations of coordinates in ${\mathbb {R}}^m$ mapping $P$ onto itself
consists of all permutations of coordinates that map facets of $P$ onto its facets. Hence, if all the facets of an $\oapl$ are known then $\gloapl$  can be calculated explicitly.
We calculated $\Pi(OAP(k,s,t,$ $\lambda))$ explicitly for all the $k,s,t$ combinations in which Theorem \ref{conj:facet} was verified. This was done by first calculating
\begin{displaymath}
 G_{k,s,t}=\{\pi|\ \mbox{there exist $\sigma $ s.t. } 
\A_{k,s,t}(\pi,\sigma)=\A_{k,s,t}\}
\end{displaymath}
where $\A_{k,s,t}$ is the constraint matrix of inequalities (\ref{eqn:Rosin})
in Theorem \ref{thm:orthogonalarray} and
$\A_{k,s,t}(\pi,\sigma)$ is the resulting
matrix when the rows of $\A_{k,s,t}$ are permuted
 according to $\sigma$ and columns according to $\pi$.
Then $H_{k,s,t}\leq \gloapl \leq G_{k,s,t}$ as $\gloapl$ must preserve the constraint matrix.
We chose calculating $G_{k,s,t}$ over directly calculating $\gloapl$ for the sake of convenience. Finding $G_{k,s,t}$ only, proved to be sufficient in all the cases we considered.
 $G_{k,s,t}$ was calculated as described in~\cite{Margot2010}, by first mapping the matrix
\[
 \left[ \begin{array}{cc}
 \0& \A_{k,s,t} \\
 \A_{k,s,t}^T &\0
\end{array} \right]
\]
to an edge colored graph and then finding its automorphism group. 
Nauty software~\cite{McKay2013}  was used to calculate the automorphism groups. In all the $k,s,t$ cases studied, it was found that 
$|G_{k,s,t}|=((s-1)!)^kk!=|H_{k,s,t}|$ implying $H_{k,s,t}=\gloapl$.
We will prove this observation after proving two lemmas.
\begin{lem}\label{lem:isometry}
Let $GG$ be the group of maps $\phi$  from  $[s]^k$   to $[s]^k$ that preserve the Hamming distance, i.e. $d(\x,\y)=d(\phi(\x),\phi(\y))$ for all $\phi \in GG$. Then $GG\cong S_s \wr S_k$. Furthermore, if $G_1$ is the subgroup of $GG$ such that 
$\tau(\1)=\1$ for all $\tau$ in $G_1$  then  $G_1\cong S_{(s-1)} \wr S_k$.
\end{lem}
\begin{pf}
Replace $\mathbb {F}_q$ with  $\mathbb{Z}_s$, the ring of integers mod $s$ and $\mathbb {F}^*_q$ with  $\mathbb{Z}_s-\{0\}$ in Theorem 5 and Theorem 6 as well as in their proofs in~\cite{Fripertinger}. Also, replace the term ``vector space" with  ``$\mathbb{Z}_s^n$". Then the resulting theorems and their proofs are still valid as the proofs never use the multiplicative invertibility of non-zero elements.
Now, replace $\{0,1,\ldots,s-1\}$ with $\{1,\ldots,s\}$ 
to get \[GG=S_{\{1,\ldots,s\}}\wr S_k\cong S_s \wr S_k\]  and \[G_1=S_{\{2,\ldots,s\}}\wr S_k\cong S_{(s-1)} \wr S_k.\]
\end{pf}


\begin{lem}\label{lem:distinctac}
If $k-t\geq 2$ then the elements of $\{a_c\}_{c=0}^t$ in Theorem \ref{thm:orthogonalarray} are all non-zero and distinct. 
\end{lem}
\begin{pf}
$\{a_c\}_{c=0}^t$ in Lemma \ref{lem:Ros} are the same $\{a_c\}_{c=0}^t$ in Theorem \ref{thm:orthogonalarray}.
$N_{\x}=\lambda/s^{k-t}$ and $N_{\y}=\lambda/s^{k-t}$ solves  the  Lemma \ref{lem:Ros} system of constraints.
Plugging in $N_{\x}=\lambda/s^{k-t}$, $N_{\y}=\lambda/s^{k-t}$ and multiplying both  sides of equations in  (\ref{eqn:Ros}) by $s^{k-t}/\lambda$ we get
\begin{equation}\label{eqn:ac}
\frac{s^{k-t}}{\lambda}a_{t-d(\z,\x)} =1+(-1)^{t-d(\z,\x)}\sum_{\y\in J_{\x} \atop d(\z,\y)>t} { d(\z,\y)-d(\z,\x)-1 \choose t-d(\z,\x)}.
\end{equation}
Taking $\z=\1$ and $\x=(2(\1_r)^T,(\1_{k-r})^T)$
equation (\ref{eqn:ac}) implies that 
\begin{equation}\label{eqn:acn}
\frac{s^{k-t}}{\lambda}a_{t-r}=1+(-1)^{t-r}\sum_{i=1}^{k-t} (s-1)^{(t-r+i)}{ k-r \choose t+i-r}
{ t+i-r-1 \choose t-r}.
\end{equation}
There are $k-t \geq 2$ positive integers inside the summation in equation  (\ref{eqn:acn}). This implies that  $a_c \neq 0$ for $0 \leq c\leq t$, and $a_ca_{c+1}<0 \Rightarrow
 a_c \neq a_{c+1}$ for  $0 \leq c\leq t-1$. Furthermore, if $a_{t-r_1} = a_{t-r_2}$ for some $1 \leq r_2<r_1\leq t$ then we must have $r_1 \equiv r_2$  (mod $2$). 
This further implies that 
{\tiny
\begin{align*}
&\sum_{i=1}^{k-t} (s-1)^{(t-r_1+i)}{ k-r_1 \choose t+i-r_1}
{ t+i-r_1-1 \choose t-r_1}=
\sum_{i=1}^{k-t}(s-1)^{(t-r_2+i)}{ k-r_2 \choose t+i-r_2}
{ t+i-r_2-1 \choose t-r_2}
\end{align*}}
for some $0 \leq r_2<r_1\leq t$. However, this is impossible as 
{\tiny
\begin{align*}
&0<(s-1)^{(t-r_1+i)}{ k-r_1 \choose t+i-r_1}
{ t+i-r_1-1 \choose t-r_1}<
(s-1)^{(t-r_2+i)}{ k-r_2 \choose t+i-r_2}
{ t+i-r_2-1 \choose t-r_2}.
\end{align*}}
\hspace {- .1in} Hence $a_{t-r_1} \neq a_{t-r_2}$ and $|a_{t-r_1}| < |a_{t-r_2}|$ for all $1 \leq r_2<r_1\leq t$.
Finally, each $a_c$ is divisible by $\lambda=a_0$ by the nature of the difference equation defining $a_c$, and $|a_c|$ is strictly increasing with $c$ as $c$ goes from $1$ to $t$. 
Hence $a_c \neq a_0=\lambda$ for $c \neq 0$.
\end{pf}
Equation (\ref{eqn:acn}) provides a closed form formula for the solution of the $a_0=\lambda$ special case of the inhomogeneous recurrence relation of degree $k-t+1$ in Lemma \ref{lem:Ros}. This recurrence relation is equivalent to a homogeneous  recurrence relation of degree $k-t+2$. Solving such an equation requires finding all complex roots of
a degree  $k-t+2$ polynomial.
  Coming up with this closed form formula for arbitrary values of $k-t$ without relating $a_c$ to the orthogonal array problem appears to be difficult.
\begin{thm}\label{thm:conj}
$\gloapl = H_{k,s,t}$.
\end{thm}
\begin{pf}
By Remark \ref{rem:remark} it suffices to consider the case $1\leq t \leq k-1$ and ($k>t+1$ or $s>2$).

Let $B_{>t}(\1):=\{\y \in [s]^k:d(\1,\y)>t\}$,   $B_{\leq t}(\1):=\{\x \in [s]^k:d(\1,\x)\leq t\}$, and $B_{t}(\1):=\{\x \in [s]^k:d(\1,\x)=t\}$.
Let $\sigma \in \gloapl$, i.e. $\sigma$ is a permutation of $B_{>t}(\1)$ which permutes the inequalities (\ref{eqn:Rosin}) of Theorem \ref{thm:orthogonalarray} by acting on the variables: $N_{\y} \mapsto N_{\sigma(\y)}$. Extend $\sigma$ to a permutation of all $[s]^k$
by assigning $\sigma(\x)=\tilde{\x}$ for $\x \in B_{\leq t}(\1)$, where $\sigma$ sends the inequality
 corresponding to $\x$ to the inequality corresponding to $\tilde{\x}$. 
 By the distinctness of the elements of $\{a_c\}_{c=0}^t$, where $c=t-d(\1,\x)$, we have
\begin{equation}\label{inq:trian}
d(\1,\x)=d(\1,\sigma(\x)) \quad \forall  \x \in B_{\leq t}(\1).
\end{equation}
For these $\x$ we also have 
\begin{equation}\label{eqn:trianph}
\sigma(J_{\x}) \cap B_{>t}(\1) = J_{\sigma(\x)} \cap B_{>t}(\1) .
\end{equation}

By (\ref{inq:trian}), $\sigma$ must map  the inequality corresponding to $\x=\1$ to itself. Thus $\sigma$ preserves the coefficients $d(\1,\y)-d(\1,\1)-1 \choose t-d(\1,\1)$ of $N_{\y}$ for $\{\y: d(\1,\y)>t\}$. Then  we must have 
\begin{equation*}\label{inq:jjjj}
d(\1,\y)=d(\1,\sigma(\y)) \quad  \forall  \y \in  B_{> t}(\1) .
\end{equation*}

For $\x \in B_{1}(\1)$, with $x_i=a\neq 1$, write $\omega(i,a)=(j,b)$ where $\sigma(\x)_j=b \neq 1$. If  $\omega(i,a)=(j,b)$
and  $\omega(\tilde{i},\tilde{a})=(j,\tilde{b})$ with $i \neq \tilde{i}$ (sharing the first component $j$), then consider an element 
$\y \in B_{>t}(\1)$ such that $y_i=a$ and  $y_{\tilde{i}}=\tilde{a}$; by equation (\ref{eqn:trianph}) above we have $\sigma(\y)_j=b = \tilde{b}$,
and $\omega$ is not injective, contradicting the injectivity of $\sigma$. Hence $i=\tilde{i}$. It follows that in the equation $ \omega(i,a)=(j,b)$, $j$ only depends on $i$ and not on $a$.

Finally, for an arbitrary $\y \in B_{>t}(\1)$, $\sigma(\y)$ is determined by the function $\omega$, since by equation (\ref{eqn:trianph})  we have $\sigma(\y)_j=b$, 
where $y_i=a \neq 1$ and $\omega(i,a)=(j,b)$;  this determines $\sigma(\y)$ uniquely, as we know that $d(\sigma(\y),\1)=d(\y,\1)$ 
and $j \neq\tilde{j}$ whenever $i \neq\tilde{i}$. Thus $\sigma$ is a permutation of columns and of non-one elements in each column, that is, $\sigma$ 
is an element  of $S_{s-1}\wr S_k$. 
\end{pf}

\begin{thm}\label{thm:marconj}
Unless  $t=0$ or $t=k$,   $G_{s,k}\cong S_{s} \wr  S_k $  in Theorem 9 of Bulutoglu and Margot~\cite{Bulutoglu2008}  is the largest subgroup of $S_{s^k}$ that sends  equations  (\ref{eqn:mar}) to themselves.
For $t=0$ or $t=k$ the largest such group is $S_{s^k}$.
\end{thm}

\begin{pf}
First, consider the case when $0<k-t<k$. Let $G$ be  the group of all coordinate permutations of $[s]^k$ which permute the rows of the constraint matrix $\A$ pertaining to equations (\ref{eqn:mar}) in the full space $\rr^{s^k}$. It was shown in Bulutoglu and Margot~\cite{Bulutoglu2008} that $G$ contains $G_{s,k}$ as a subgroup. Hence, it suffices to show that $G \leq G_{s,k}$. First, we show that every element in $G$ is an isometry. 

Let $\sigma \in G$ and $\z \in [s]^k$.
Since  $G_{s,k}$ acts transitively on  $[s]^k$ and itself consists of isometries, we may assume $\sigma(\z)=\z$. 
Consider the set $S$ of all rows of $\A$ in which $\z$ appears. Since $\sigma(\z)=\z$, $\sigma$ must permute $S$.
Let $\w\in [s]^k$ be such that $d(\z,\w)=i$, where  $i \in \{0,1,\ldots,k-t\}$. Then $\w$ appears in exactly ${k-i\choose t}$ rows in $S$.
Since $\sigma$ preserves $S$, $\sigma(\w)$ also appears in exactly ${k-i\choose t}$ rows in $S$. Then we must have  $d(\z,\sigma(\w))=i$ as no element is repeated in the set
$\{{k-i\choose t}\}_{i=0}^{k-t}$.
Hence $d(\z,\w)=d(\sigma(\z),\sigma(\w))$ for all $\z$ and $\w$ such that $0\leq d(\z,\w)\leq k-t$.


Now, for any $\z$ and  $\y$ such that $d(\z,\y)=k-t+1$ there exists $\y_1$ such that $d(\z,\y)=d(\z,\y_1)+d(\y_1,\y)$, where
$d(\z,\y_1)=k-t$ and $d(\y_1,\y)=1 \leq k-t$.
 By the triangle inequality $d(\sigma(\z),\sigma(\y)) \leq d(\sigma(\z),\sigma(\y_1)) +d(\sigma(\y_1),\sigma(\y))=k-t+1=d(\z,\y)$.
By repeating the same argument for $\z$ and $\y$ such that $d(\z,\y)=k-t+i$ for $i=2,\ldots,t$ we get
\begin{equation}\label{inq:trianph}
d(\sigma(\z),\sigma(\y)) \leq d(\z,\y)
\end{equation}
 for all $\z,\y$ and $\sigma$. Let $\z'=\sigma(\z)$, $\y'=\sigma(\y)$ and $h=\sigma^{-1}$ then
\begin{equation}\label{inq:trianh}
d(\z',\y') \leq d(h(\z'),h(\y'))
\end{equation}
 for all  $\z',\y'$ and $h$. Combining inequalities (\ref{inq:trianph}) and (\ref{inq:trianh}) we get
\begin{equation*}\label{inq:trianphfin}
d(\sigma(\z),\sigma(\y)) = d(\z,\y)
\end{equation*}
for all $\z,\y \in [s]^k$ and $\sigma \in  G$. Hence, by Lemma \ref{lem:isometry}, an isomorphic copy of $G$ is contained in 
$ S_{s}\wr S_{k}.$
 Now, this implies that $G \leq G_{s,k} \cong S_{s} \wr  S_{k} $.
The cases $t=0$ and $t=k$ are easy to see.
\end{pf}

\begin{center} {\sc Acknowledgements} \end{center}

The views expressed in this article are
those of the authors and do not reflect the official policy or
position of the United States Air Force, Department of Defense, or
the US Government.

This research was supported by the AFOSR grant F1ATA03039J001.

\bibliographystyle{elsarticle-num}
\bibliography{BibliographyDORos}

\newpage
\tiny
\label{tabel:long}
\begin{longtable}{ccccc}
\caption{Formulation Comparisons}\\
\hline
\multicolumn{1}{c}{}&\multicolumn{1}{c}{Bulutoglu and}&
\multicolumn{1}{c}{After Deleting }&\multicolumn{1}{c}{Bulutoglu and}&
\multicolumn{1}{c}{After Deleting }\\
\multicolumn{1}{c}{}&\multicolumn{1}{c}{Margot [2] }&
\multicolumn{1}{c}{Variables }&\multicolumn{1}{c}{Margot [2]}&
\multicolumn{1}{c}{Variables }\\
\multicolumn{1}{c}{}&\multicolumn{1}{c}{Number of}&
\multicolumn{1}{c}{ Number of}&\multicolumn{1}{c}{  Formulation}&
\multicolumn{1}{c}{  Formulation}\\
\multicolumn{1}{c}{OA($N,k,s,t$)}&\multicolumn{1}{c}{$\oan$s}&
\multicolumn{1}{c}{$\oan$s}&\multicolumn{1}{c}{Times (sec.)}&
\multicolumn{1}{c}{ Times (sec.)}\\
\hline
\endfirsthead

\hline \multicolumn{5}{c}%
{\tablename\ \thetable{} -- {\rm continued \ from \ previous \ page}}\\
\hline
\multicolumn{1}{c}{}&\multicolumn{1}{c}{Bulutoglu and}&
\multicolumn{1}{c}{After Deleting }&\multicolumn{1}{c}{Bulutoglu and}&
\multicolumn{1}{c}{After Deleting }\\
\multicolumn{1}{c}{}&\multicolumn{1}{c}{Margot [2] }&
\multicolumn{1}{c}{Variables }&\multicolumn{1}{c}{Margot [2]}&
\multicolumn{1}{c}{Variables }\\
\multicolumn{1}{c}{}&\multicolumn{1}{c}{Number of}&
\multicolumn{1}{c}{ Number of}&\multicolumn{1}{c}{  Formulation}&
\multicolumn{1}{c}{  Formulation}\\
\multicolumn{1}{c}{OA($N,k,s,t$)}&\multicolumn{1}{c}{$\oan$s}&
\multicolumn{1}{c}{$\oan$s}&\multicolumn{1}{c}{Times (sec.)}&
\multicolumn{1}{c}{ Times (sec.)}\\
\hline
\endhead
\hline \multicolumn{5}{r}{{\rm Continued \ on \ next \ page}}\\
\hline
\endfoot
 \hline
\endlastfoot
\hline
OA(20,6,2,2) &75  & 3069 & 1.42 & 63.99 \\
OA(20,7,2,2) &474 & 51695 & 13.4 &  2578.82 \\
OA(20,8,2,2) &1603 & 383729 & 108.96 &  66377 \\
OA(20,9,2,2) &2477 & 1157955 & 484.55 &  879382 \\
OA(20,10,2,2) &2389 & $\geq28195$ & 1683.95 & $\geq 37214 $\\
OA(24,5,2,2) &63  & 723 & 1.07 &  18.36 \\
OA(24,6,2,2) &1350  & 62043 & 22.03 &  1381.39 \\
OA(24,7,2,2) &57389 & 6894001 & 1720.96 &  428220 \\
OA(24,8,2,2) &1470157 & 4505018 & 99738 &  653671 \\
OA(24,9,2,2) &3815882 & -- & 763643 & -- \\
OA(24,5,2,3) &1  & 2 & 0.13 &  11.64 \\
OA(24,6,2,3) &2  & 5 & 0.25 &  11.67 \\
OA(24,7,2,3) &1  & 5 & 0.32 &  16.04 \\
OA(24,8,2,3) &1  & 6 & 1 &  22.88 \\
OA(24,9,2,3) &1  & 6 & 5.9 &  44.02 \\
OA(24,10,2,3) &1 & 5 & 55.49&  128.95 \\
OA(24,11,2,3) &1 & 3 & 519.62&  460.59 \\
OA(32,6,2,3) &10 & 31 & 1.85 &  12.2 \\
OA(32,7,2,3) &17 & 76 & 1.82 &  16.13 \\
OA(32,8,2,3) &33  & 194 & 6.59 &  77.49 \\
OA(32,9,2,3) &34  & 364 & 23.75 &  658.38 \\
OA(32,10,2,3) &32  & 561 & 102.39 &  7338 \\
OA(32,11,2,3) &22  &$\geq 441$ & 560.29 & $\geq36463$ \\
OA(40,6,2,3) &9  & 65 & 0.52 &  12.92 \\
OA(40,7,2,3) &25 & 580 & 2.01 &  40.68 \\
OA(40,8,2,3) &105& 6943 & 19.71 & 4178 \\
OA(40,9,2,3) &213 & 43713 & 206.25 &  260919 \\
OA(40,10,2,3) &353 & $\geq1511$ & 1764.73 & $ \geq36279$ \\
OA(48,6,2,3) &45 & 355 & 2.01 &  18.27 \\
OA(48,7,2,3) &397 & 13469 & 33.73 & 862.1 \\
OA(48,8,2,3) &8383 & 896963 & 2231.77 & 552154 \\
OA(54,5,3,3) &4  & 49 & 1.9 & 36.01 \\
OA(54,6,3,3) &0  & 0 & 17.14 &  167.07 \\
OA(56,6,2,3) &86  & 1393 & 4.44 &  36.02 \\
OA(56,7,2,3) &4049  & 285184 & 443.4 & 20415 \\
OA(64,7,2,4) &7  & 21 & 98.83 & 15.45 \\
OA(64,8,2,4) &3  & 10 & 12.17 & 23.39 \\
OA(80,6,2,4) &1  & 6 & 0.52 & 11.86 \\
OA(80,7,2,4) &0  & 0 & 0.37 & 15.01 \\
OA(81,5,3,4) &1  & 2 & 15.75 &  19.56 \\
OA(96,7,2,4) &4  & 31 & 3.14 &  15.41 \\
OA(96,8,2,4) &0  & 0 & 2.28 &  60.39 \\
OA(112,6,2,4) &3  & 25 & 1.24 &  12.7 \\
OA(112,7,2,4) &0  & 0 & 1.24 &  17.36 \\
OA(144,8,2,4) &20 & 3392 & 1792.82 & 1535314 \\
OA(162,6,3,4) &0  & 0 & 19.8 & 266.8 \\
\hline
\end{longtable}
\end{document}